\documentclass[11pt]{amsart}
\usepackage{graphicx}
\usepackage{amssymb, amsmath}
\newtheorem{theorem}{Theorem}[section]
\newtheorem{corollary}[theorem]{Corollary}
\newtheorem{lemma}[theorem]{Lemma}

\theoremstyle{definition}
\newtheorem{definition}[theorem]{Definition}

\theoremstyle{remark}

\begin{document}
\title{On the Automorphisms and Representations of Polyadic Groups}
\author{H. Khodabandeh}
\address{Department of Pure Mathematics,  Faculty of Mathematical
Sciences, University of Tabriz, Tabriz, Iran}
\author{\sc M. Shahryari}
\thanks{{\scriptsize
\hskip -0.4 true cm MSC(2010): 20N15
}}

\address{ Department of Pure Mathematics,  Faculty of Mathematical
Sciences, University of Tabriz, Tabriz, Iran }
\email{mshahryari@tabrizu.ac.ir}
\date{\today}

\begin{abstract}
Using a unified method, we determine the structure of automorphisms and representations of arbitrary polyadic groups. More precisely, for a polyadic group $(G, f)=der_{\theta, b}(G, \cdot)$, we obtain a complete description of automorphisms and representations of $(G,f)$ in terms of automorphisms and representations of the binary group $(G, \cdot)$.
\end{abstract}

\maketitle

\section{Introduction}
A non-empty set $G$ together with an $n$-ary operation $f:G^n\to G$
is called an {\it $n$-ary quasi-group} or a {\it polyadic quasi-group}, if for all
$x_0,x_1,\ldots,x_n\in G$ and for any fixed $i\in\{1,\ldots,n\}$, there
exists a unique  element $y\in G$, such that
$$
f(x_1, \ldots, x_{i-1},y,x_{i+1}, \ldots, x_n)=x_0.
$$
In the binary case (i.e., for $n=2$), a polyadic quasi-group is just usual
quasi-group. If the operation $f$ is also associative, then we call $(G,f)$ an {\it n-ary group} or a {\it polyadic
group}.

According to a general convention used in the theory of $n$-ary systems,
the sequence of elements $x_i,x_{i+1},\ldots,x_j$ is denoted by
$x_i^j$. In the case $j<i$ it is the empty symbol. If
$x_{i+1}=x_{i+2}=\ldots=x_{i+t}=x$, then instead of
$x_{i+1}^{i+t}$ we write $\stackrel{(t)}{x}$. In this convention
$f(x_1,\ldots,x_n)= f(x_1^n)$ and
$$
 f(x_1,\ldots,x_i,\underbrace{x,\ldots,x}_{t},x_{i+t+1},\ldots,x_n)=
 f(x_1^i,\stackrel{(t)}{x},x_{i+t+1}^n) .
$$

As an example of a polyadic quasi-group, suppose $(G, \cdot)$ is an ordinary quasi-group (so the equations
$ax=b$ and $xa=b$ have unique solutions for $x$). Let $\alpha_1, \ldots, \alpha_n$  be arbitrary automorphisms of
$(G, \cdot)$ and $b\in G$ be fixed. If we define
$$
f(x_1^n)=\alpha_1(x_1)\ldots\alpha_n(x_n) b,
$$
then $(G, f)$ becomes an $n$-ary quasi-group. A polyadic quasi-group of this type, is called a {\it linear polyadic quasi-group}.
As an special case, if $(G, \cdot)$ is an ordinary group, $\alpha_1= \cdots= \alpha_n=identity$ and $b=1$, then we have
$$
 f(x_1^n)=x_1x_2 \ldots x_n,
$$
and $(G, f)$ is an $n$-ary group which is called the $n$-ary group, derived from $(G, \cdot)$ and it is denoted by $der^n(G, \cdot)$. More generally,
suppose $(G, \cdot)$ is an ordinary group and $\theta$ is an automorphism with a fixed point $b\in G$, such that $\theta^{n-1}(x)=bxb^{-1}$ for all $x\in G$. Now, define
$$
f(x_1^n)=x_1\theta(x_2)\theta^2(x_3)\ldots \theta^{n-1}(x_n)b.
$$
It can be shown that $(G,f)$ is an $n$-ary group which we denote it by $der_{\theta, b}(G, \cdot)$ and we call it $(\theta, b)$-{\it derived polyadic
 group} from $(G, \cdot)$. It is proved that the converse is also true, i. e. every polyadic 
 group can be uniquely expressed as  $(G, f)=der_{\theta, b}(G, \cdot)$. This is the content of Hossz\'u-Gluskin Theorem, 
 which is formulated as follows.

\begin{theorem}
Let $(G,f)$ be an $n$-ary group. Then

\begin{enumerate}
\item[$(1)$] \ on $G$ one can define an operation $\cdot$ such
that $(G,\cdot)$ is a group,
\item[$(2)$] \ there exist an automorphism $\theta$ of $(G,\cdot)$ and $b\in G$, such
that $\theta(b)=b$,
\item[$(3)$] \ $\theta^{n-1}(x)=b x b^{-1}$, for every $x\in
G$,
\item[$(4)$] \ $
f(x_1^n)=x_1\theta(x_2)\theta^2(x_3)\cdots\theta^{n-1}(x_n)b$, for all $x_1,\ldots,x_n\in G$.
\end{enumerate}
\end{theorem}
For a proof, see  \cite{DG} or \cite{DM}.

From the definition of an $n$-ary group $(G,f)$ we can directly see
that for every $x\in G$, there exists only one $y\in G$ satisfying
the equation
$$
f(\stackrel{(n-1)}{x},y)=x .
$$
This element is called {\it skew} to $x$ and is denoted by
$\overline{x}$. In a ternary group ($n=3$) derived from a binary
group $(G,\cdot)$, the skew element coincides with the inverse
element in $(G,\cdot)$. Thus, in some sense, the skew element is a
generalization of the inverse element in binary groups. Nevertheless, the concept of skew elements plays a crucial role in
the theory of $n$-ary groups. Namely, as D\"ornte proved, the following theorem (see \cite{Dor}).

\begin{theorem}
In any $n$-ary group $(G,f)$ the following identities
$$
f(\stackrel{(i-2)}{x},\overline{x},\stackrel{(n-i)}{x},y)=
f(y,\stackrel{(n-j)}{x},\overline{x},\stackrel{(j-2)}{x})=y,
$$
$$
f(\stackrel{(k-1)}{x},\overline{x},\stackrel{(n-k)}{x})=x
$$
 hold for all $\,x,y\in G$, $\,2\leqslant i,j\leqslant
n$ and $1\leqslant k\leqslant n$.
\end{theorem}

Fixing in an $n$-ary operation $f$ all inner elements
$a_2,\ldots,a_{n-1}$ we obtain a new binary operation
$$
x\cdot y=f(x,a_2^{n-1},y).
$$
Such obtained groups $(G,\cdot)$ is called a {\it retract} of
$(G,f)$. Choosing different elements $a_1,\ldots,a_{n-1}$ we
obtain different retracts.
Retracts of a fixed $n$-ary group are isomorphic (cf.
\cite{DM}). So, we can consider only retracts of the form
$$
x\cdot y=f(x,\stackrel{(n-2)}{a},y).
$$
Such retracts will be denoted by $Ret_a(G,f)$. The identity of the
group $Ret_a(G,f)$ is $\overline{a}$. One can verify that the
inverse element to $x$ has the form
\begin{equation*}
x^{-1}=f(\overline{a},\stackrel{(n-3)}{x},\overline{x},\overline{a}).
\end{equation*}
The binary group, which we said about in the theorem 1.1, is in fact $Ret_a(G, f)$ and the automorphism $\theta$ is defined as
$\theta(x)=(\bar{a}, x, \stackrel{(n-2)}{a})$.

Binary retracts of an $n$-ary group $(G,f)$ are commutative only
in the case when there exists an element $a\in G$ such that
$$
f(x,\stackrel{(n-2)}{a},y)=f(y,\stackrel{(n-2)}{a},x)
$$
holds for all $x,y\in G$. An $n$-ary group with this property is
called {\em semiabelian}. It satisfies the identity
\begin{equation*}
 f(x_1^n)=f(x_n,x_2^{n-1},x_1)
\end{equation*}
(see \cite{Rem}).

One can prove (cf. \cite{GG}) that a semiabelian $n$-ary group
is {\it medial}, i.e., it satisfies the identity
$$
f(f(x_{11}^{1n}),f(x_{21}^{2n}),\ldots,f(x_{n1}^{nn}))=f(f(x_{11}^{n1}),f(x_{12}^{n2}),\ldots,
f(x_{1n}^{nn})).
$$
In such $n$-ary groups
\begin{equation*}
\overline{f(x_1^n)}=f(\overline{x}_1,\overline{x}_2,\overline{x}_3,\ldots,\overline{x}_n)
\end{equation*}
for all $x_1,\ldots,x_n\in G$.

The idea of investigations of such polyadic group  goes  back to E.
Kasner's lecture \cite{Kas} at the fifty-third annual meeting of the
American Association for the Advancement of Science in 1904. But the
first paper concerning the theory of $n$-ary groups was written
(under inspiration of Emmy Noether) by W. D\"ornte in 1928 (see
\cite{Dor}). In 1940 E. L. Post   published
an extensive study of n-groups in which the well-known Post's Coset Theorem
appeared, see \cite{Post}.

Representation theory of polyadic groups is investigated by W. Dudek and M. Shahryari in \cite{Shah}.
In \cite{Shah2}, the second author,  studied characters of finite polyadic groups. In this article, using a unified method (homotopies and autotopies of polyadic groups), we study the structure of automorphisms and representations of polyadic group. More precisely, if $(G, f)=der_{\theta, b}(G, \cdot)$, we express an automorphism (a representation) of $(G, f)$, as a product of a certain automorphism (representation) of the binary group
$(G, \cdot)$ and a translation.

\section{Automorphisms}
Suppose $(G,f)$ and $(H, h)$ are two $n$-ary quasi-groups. Suppose there are maps
$\alpha_1, \ldots, \alpha_{n+1}$ from $G$ to $H$, such that
$$
\alpha_{n+1}(f(x_1, \ldots, x_n))=h(\alpha_1(x_1), \ldots, \alpha_n(x_n)).
$$
Then we say that $T=(\alpha_1, \ldots, \alpha_{n+1})$ is a {\it homotopy} from $(G, f)$ to $(H, h)$. If all maps, $\alpha_i$, are bijections,
then $T$ is an {\it isotopy}. We call $T$ an {\it autotopy} of $(G, f)$, if $(G,f)=(H, h)$ and $T$ is an isotopy.
The set of all autotopies of an $n$-ary
quasi-group $(G, f)$ is denoted by
$\mathfrak{T}(G, f)$. If we define a binary operation
$$
(\alpha_1, \ldots, \alpha_{n+1})\circ (\beta_1, \ldots, \beta_{n+1})=(\alpha_1\beta_1, \ldots, \alpha_{n+1}\beta_{n+1}),
$$
then $\mathfrak{T}(G, f)$ becomes a group. In general, if $T=(\alpha_1, \ldots, \alpha_{n+1})$ is a  homotopy from $(G, f)$ to $(H, h)$ and
$S=(\beta_1, \ldots, \beta_{n+1})$ is another homotopy from $(K, g)$ to $(H, h)$, then we define their composition as
$$
(\alpha_1, \ldots, \alpha_{n+1})\circ (\beta_1, \ldots, \beta_{n+1})=(\alpha_1\beta_1, \ldots, \alpha_{n+1}\beta_{n+1}).
$$
It is easy to see that if $T=(\alpha_1, \ldots, \alpha_{n+1})$ is an isotopy from $(G, f)$ to $(H, h)$, then we have
$$
\mathfrak{T}(G, f)=T^{-1}\circ \mathfrak{T}(H, h)\circ T.
$$
A map $\psi: G\to H$ is a {\it homomorphism} from $(G, f)$ to $(H, h)$, if $T=(\psi, \ldots, \psi)$ is a homotopy. If
$$
T=(\psi, \ldots, \psi)\in  \mathfrak{T}(G, f),
$$
then we say that $\psi$ is an automorphism of $(G,f)$. For example, let $(G, f)=der_{\theta, b}(G, \cdot)$. Then for any $a\in G$, the map
$\psi(x)=a\theta(xa^{-1})b^{-1}$ is an automorphism of $(G, f)$. The group of automorphisms of an arbitrary $(G, f)$, will be denoted by
$Aut(G, f)$.

Let $(G, \cdot)$ be an ordinary quasi-group and $a\in G$. We define the {\it left} and the {\it right translation maps}, $L_a$ and $R_a$ by
$L_a(x)=ax$ and $R_a(x)=xa$, respectively. If $G$ is a group, we denote by $I_a$, the inner automorphism $I_a=(x)a^{-1}xa$. For the proof of
the next theorem, see \cite{Mar}.

\begin{theorem}
Suppose $(G, \cdot)$ is an ordinary group and $(G, f)=der^n(G, \cdot)$. Then $T$ is an autotopy of
$(G, f)$, if and only if
$$
T=(L_{a_1}I_{a_1}, L_{a_2}I_{a_1a_2}, \ldots, L_{a_n}I_{a_1\ldots a_n}, R_{a_1\ldots a_n})\circ (\phi, \ldots, \phi),
$$
where $a_1^n\in G$ and $\phi\in Aut(G, \cdot)$. Further, the above decomposition is unique.
\end{theorem}

Now, we are going to determine the structure of automorphism of a polyadic group of type $(G, f)=der_{\theta}(G, \cdot)$, (i.e., $b=1$). The general case, will be proved later. In what follows, $\varepsilon$ denotes the identity map.

\begin{theorem}
Let $(G, f)=der_{\theta}(G, \cdot)$, where $(G, \cdot)$ is an ordinary group and $\theta$ is an ordinary automorphism with $\theta^{n-1}=\varepsilon$.
Then
$$
Aut(G, f)=\{ R_a\phi: \bar{a}=a, \phi\in Aut(G, \cdot), [\theta, \phi]=I_a\},
$$
where the bracket denotes the ordinary commutator $\theta\phi\theta^{-1}\phi^{-1}$.
\end{theorem}

{\bf Proof}. Define a second operation $g$ on $G$ by $g(x_1^n)=x_1x_2\ldots x_n$. Now, $T=(\varepsilon, \theta, \theta^2, \ldots, \theta^{n-2}, \varepsilon, \varepsilon)$ is an isotopy between $(G, f)$ and $(G, g)$. So we have
$$
\mathfrak{T}(G, f)=T^{-1}\circ \mathfrak{T}(G, g)\circ T.
$$
Using 2.1, we have also,
\begin{eqnarray*}
\mathfrak{T}(G, g)&=&\{ (L_{a_1}I_{a_1}, L_{a_2}I_{a_1a_2}, \ldots, L_{a_n}I_{a_1\ldots a_n}, R_{a_1\ldots a_n})\circ (\phi, \ldots, \phi)\\
                  & &\ \ \ \ \ \ \ \ \ \ \ \ \ \ \ \ \ \ \ \ \ \ \ \ \ \ \ \ \ \ \ \ \ \ \ \ \ \ \ \ : a_1^n\in G, \phi\in Aut(G, \cdot)\}.
\end{eqnarray*}
Hence, every element of $\mathfrak{T}(G, f)$ can be written uniquely as
$$
(L_{a_1}I_{a_1}\phi, \theta^{-1}L_{a_2}I_{a_1a_2}\phi\theta, \ldots,
 \theta^{1-n}L_{a_n}I_{a_1\ldots a_n}\phi\theta^{n-1},R_{a_1\ldots a_n}\phi),
$$
where $a_1^n\in G$ and $\phi\in Aut(G, \cdot)$. It is clear that $Aut(G, f)$ is precisely the set of all $\psi$'s such that $(\psi, \ldots, \psi)\in
\mathfrak{T}(G, f)$. So, we must determine all elements of $\mathfrak{T}(G, f)$, with the equal entries. Therefore, suppose we have
\begin{eqnarray*}
\psi&=&L_{a_1}I_{a_1}\phi\\
    &=&\theta^{-1}L_{a_2}I_{a_1a_2}\phi\theta\\
    &\vdots&\\
    &=&\theta^{2-n}L_{a_{n-1}}I_{a_1\ldots a_{n-1}}\phi\theta^{n-2}\\
    &=& L_{a_n}I_{a_1\ldots a_n}\phi\\
    &=&R_{a_1\ldots a_n}\phi.
\end{eqnarray*}
Suppose also $a=a_1$. Then, clearly $\psi=R_a\phi$, thus we prove that $\bar{a}=a$ and $[\theta, \phi]=I_a$. We, have the following steps.

i- The equality of the first and the last entry ($R_a\phi=R_{aa_2\ldots a_n}\phi$), implies $a_2a_3\cdots a_n=1$.

ii- The equality of the first and $n$-th entry ($R_a\phi=L_{a_n}I_{aa_2\ldots a_n}\phi$), implies $a_n=a$.

iii- Now, we use the equality of the second and the third entries. We have $R_a\phi=\theta^{-1}L_{a_2}I_{aa_2}\phi\theta$, so for all $x$,
\begin{eqnarray*}
\phi(x)a&=&\theta^{-1}(a_2a_2^{-1}a^{-1}\phi(\theta(x))aa_2)\\
        &=&\theta^{-1}(a^{-1})\phi(\theta(x))\theta^{-1}(a)\theta^{-1}(a_2).
\end{eqnarray*}
Hence, if we put $x=1$, then we obtain $\theta(a)=a_2$.

iv- Continuing this way, using other equalities, we conclude that for all $i$, $a_i=\theta^{i-1}(a)$. Hence
\begin{eqnarray*}
f(\stackrel{(n)}{a})&=&a\theta(a)\ldots \theta^{n-1}(a)\\
                    &=&aa_2\ldots a_n\\
                    &=&a.
\end{eqnarray*}
So $\bar{a}=a$. Now, we have $a\theta(a)\theta^2(a)\ldots \theta^{n-2}(a)=1$, hence applying the equality of the first and $(n-2)$-th entries,
we obtain
\begin{eqnarray*}
R_a\phi(x)&=&\theta^{2-n}L_{\theta^{n-2}(a)}I_{a\theta(a)\ldots \theta^{n-2}(a)}\phi\theta^{n-2}(x)\\
          &=&\theta^{2-n}L_{\theta^{n-2}(a)}\phi\theta^{n-2}(x)\\
          &=&\theta^{2-n}(\theta^{n-2}(a)\phi(\theta^{n-2}(x)))\\
          &=&a\theta^{2-n}(\phi(\theta^{n-2}(x))).
\end{eqnarray*}
So, we have $ \theta^{2-n}\phi\theta^{n-2}(x)=a^{-1}\phi(x)a$. But, $\theta^{n-2}=\theta^{-1}$, hence we have
$$
\theta \phi \theta (x)=a^{-1}\phi(x)a,
$$
which is equivalent to $[\theta, \phi]=I_a$.

Conversely, suppose $a\in G$ has the property $\bar{a}=a$. Let $\phi\in Aut(G, \cdot)$ such that $[\theta, \phi]=I_a$. We prove that
$R_a\phi$ is an automorphism of $(G, f)$. Since $\bar{a}=a$, we know that the following $(n+1)$-tuple is an autotopy of $(G, f)$;
$$
(R_a\phi, \theta^{-1}L_{\theta(a)}I_{a\theta(a)}\phi\theta, \ldots,
 \theta^{1-n}L_{\theta^{n-1}(a)}I_{a}\phi\theta^{n-1},R_a\phi).
$$
It is enough to show that all entries of the above autotopy are equal. Using $[\theta, \phi]=I_a$, we obtain
$$
\theta^{-i}\phi\theta^i(x)=\theta^{-1}(a)\cdots \theta^{-1}(a)\phi(x)\theta^{-1}(a^{-1})\cdots \theta^{-i}(a^{-1}),
$$
for all $i$. A few computations show that the above relation is equivalent to
$$
R_a\phi=\theta^{-i}L_{\theta^i(a)}I_{a\theta(a)\cdots \theta^i(a)}\phi\theta^i.
$$
Hence, all entries of our autotopy are equal and so $R_a\phi$ is an automorphism.\\

\begin{definition}
An element $a\in G$ with the property $\bar{a}=a$ is called an {\it idempotent}. If $(G, f)=der_{\theta}(G, \cdot)$, then the set of all idempotents which are also in the center of $(G, \cdot)$, will be denoted by $Z^{\ast}(G)$. It is easy to see that $Z^{\ast}(G)$ is a subgroup of $(G, \cdot)$.
\end{definition}

\begin{corollary}
Suppose $(G, f)=der_{\theta}(G, \cdot)$. Then we have
$$
\frac{Aut(G,f)}{ Z^{\ast}(G)}\hookrightarrow Aut(G, \cdot).
$$
If further, all idempotents of $(G, f)$ are central, then we have
$$
Aut(G, f)\cong C_{Aut(G, \cdot)}(\theta)\ltimes Z^{\ast}(G).
$$
\end{corollary}

{\bf Proof}. Note that an automorphism $\psi \in Aut(G,f)$ can be uniquely expressed in the form $\psi=R_a\phi$, so we can
define a map $q:Aut(G, f)\to Aut(G, \cdot)$ by $q(R_a\phi)=\phi$. It is easy to check that $(R_a\phi)(R_{a^{\prime}}\phi^{\prime})=
R_{\phi(a^{\prime})a}\phi\phi^{\prime}$. Hence $q$ is a group homomorphism. Clearly, its kernel is $Z^{\ast}(G)$, so
$$
\frac{Aut(G,f)}{ Z^{\ast}(G)}\hookrightarrow Aut(G, \cdot).
$$
Now, let all idempotents of $(G, f)$ be central, so
$$
Aut(G, f)=\{ R_a\phi: a\in Z^{\ast}(G), \phi\in C_{Aut(G, \cdot)}(\theta)\}.
$$
Define an action of $C_{Aut(G, \cdot)}(\theta)$ on $Z^{\ast}(G)$ by the rule $\phi.a=\phi(a)$. One can check that the map
$$
\lambda:C_{Aut(G, \cdot)}(\theta)\ltimes Z^{\ast}(G)\to Aut(G, f)
$$
with the rule $\lambda (\phi, a)=R_a\phi$ is an isomorphism. This completes the proof.\\

\begin{corollary}
Suppose $(G, f)=der_{\theta}(G, \cdot)$ is a medial polyadic group. Then
$$
Aut(G, f)\cong C_{Aut(G, \cdot)}(\theta)\ltimes Z^{\ast}(G).
$$
\end{corollary}

\begin{corollary}
Let $(G, \cdot)$ be an abelian group and $(G,f)=der^n(G, f)$. Then
$$
Aut(G, f)=Aut(G, \cdot)\ltimes \bar{Z}(G),
$$
where $\bar{Z}(G)=\{ a\in G: a^{n-1}=1\}$.
\end{corollary}

We are ready now, to talk about the structure of automorphisms of a polyadic group in the general form
$(G, f)=der_{\theta, b}(G, \cdot)$. Note that $\theta$ and $b$ satisfy the conditions of the theorem 1.1. \\

\begin{theorem}
Suppose $(G, f)=der_{\theta, b}(G, \cdot)$. Then we have
$$
Aut(G, f)=\{ R_a\phi: \phi\in Aut(G, \cdot), f(\stackrel{(n)}{a})=\phi(b)a, [\theta, \phi]=I_a\},
$$
\end{theorem}

{\bf Proof}. Let $a\in G$ and $\phi\in Aut(G, \cdot)$ satisfy $f(\stackrel{(n)}{a})=\phi(b)a$ and  $[\theta, \phi]=I_a$. We show that
$\psi=R_a\phi$ is an automorphism of $(G,f)$. For any $x\in G$ and $1\leq i\leq n-1$, we have
$$
\phi\theta^i(x)=a\theta(a)\cdots \theta^{i-1}(a)\theta^i(\phi(x))\theta^{i-1}(a^{-1})\cdots\theta(a^{-1})a^{-1}.
$$
Hence
\begin{eqnarray*}
R_a\phi(f(x_1^n))&=&R_a\phi(x_1\theta(x_2)\cdots\theta^{n-1}(x_n)b)\\
                 &=&\phi(x_1)\phi(\theta(x_2))\cdots\phi(\theta^{n-1}(x_n))\phi(b)a\\
                 &=&\phi(x_1)(a\theta(\phi(x_2))a^{-1})(a\theta(a)\theta^2(a)(\phi(x_3))\theta^2(a^{-1})\theta(a^{-1})a^{-1})\cdots\\
                 & & (a\theta(a)\cdots \theta^{n-2}(a)\theta^{n-1}(\phi(x))\theta^{n-2}(a^{-1})\cdots\theta(a^{-1})a^{-1})f(\stackrel{(n)}{a})\\
                 &=&(\phi(x_1)a)(\theta(\phi(x_2))\theta(a))\cdots (\theta^{n-1}(\phi(x_n))\theta^{n-1}(a))b\\
                 &=&f(R_a\phi(x_1), \ldots, R_a\phi(x_n)).
\end{eqnarray*}
Now, to show that every automorphism of $(G,f)$ has the required form, suppose $(G,g)=der^n(G, \cdot)$. It is clear that
$$
T=(\varepsilon, \theta, \theta^2, \ldots, \theta^{n-1}, R_{b^{-1}})
$$
is an isotopy between $(G, f)$ and $(G,g)$. So we can do the same argument as in 2.2, to complete the proof.\\

\begin{lemma}
Let $(G, f)=der_{\theta, b}(G, \cdot)$ and suppose $u\in G$ is a central idempotent. Then the map
$$
R_{u}:der_{\theta}(G, \cdot)\to (G, f),
$$
is an isomorphism.
\end{lemma}

\begin{corollary}
If $(G, f)=der_{\theta, b}(G, \cdot)$ has a central idempotent, then
$$
Aut(G, f)\cong Aut(der_{\theta}(G, \cdot)).
$$
\end{corollary}

\section{Homotopy and the structure of homomorphisms }
Employing the same method as in the section 2, we want to determine the structure a homomorphism between two
polyadic groups. We apply the results of this section, to investigate representations of polyadic groups in the section 4.

\begin{lemma}
Suppose $(G, \cdot)$ and $(H, \ast)$ are two ordinary groups. Then every homomorphism
$$
\psi:der^n(G, \cdot)\to der^n(H, \ast)
$$
can be uniquely decomposed as $\psi=R_a\phi$ such that\\
i- \ \ \ $\phi:(G, \cdot)\to (H, \ast)$ is an ordinary homomorphism,\\
ii-\ \ \ $a\in C_H(\phi(G))$,\\
iii-\ \  $a^{n-1}=1$.\\
The converse is also true.
\end{lemma}

{\bf Proof}. Suppose $\psi$ is given and let $a=\psi(1)$. Then $\psi(1^n)=a$ and so $a^{n-1}=1$. Now, define
$\phi=R_a^{-1}\psi$. We have
\begin{eqnarray*}
\phi(xy)&=&\psi(xy)\ast a^{-1}\\
        &=&\psi(x\cdot 1^{n-2}\cdot y)\ast a^{-1}\\
        &=&\psi(x)\ast a^{n-2}\ast \psi(y)\ast a^{-1}\\
        &=&\psi(x)\ast a^{-1}\ast \psi(y)\ast a^{-1}\\
        &=&\phi(x)\ast \phi(y).
\end{eqnarray*}
So, $\phi$ is an ordinary homomorphism. Further
\begin{eqnarray*}
\psi(x)&=&\psi(1\cdot x\cdot 1^{n-2})\\
       &=&a\ast \psi(x)\ast a^{-1},
\end{eqnarray*}
so $a\ast \psi(x)=\psi(x)\ast a$, which implies that $a\in C_H(\phi(G))$. Conversely, suppose $\psi=R_a\phi$, such that $a$ and
$\phi$ satisfy the above three conditions. We have
\begin{eqnarray*}
\psi(x_1x_2\ldots x_n)&=& \phi(x_1x_2\ldots x_n)\ast a\\
                      &=&\phi(x_1)\ast \phi(x_2)\ast \cdots \ast \phi(x_n)\ast a^n\\
                      &=&\phi(x_1)\ast a\ast \phi(x_2)\ast a\ast \cdots \ast\phi(x_n)\ast a\\
                      &=&\psi(x_1)\ast \psi(x_2)\ast \cdots \ast \psi(x_n).
\end{eqnarray*}
This completes the proof.\\

The following theorem is a generalization of 2.1.\\

\begin{theorem}
Suppose $(G, \cdot)$ and $(H, \ast)$ are two ordinary groups. Then every homotopy
$$
der^n(G, \cdot)\to der^n(H, \ast)
$$
can be  decomposed as
$$
(L_{a_1}I_{a_1}, L_{a_2}I_{a_1\ast a_2}, \ldots, L_{a_n}I_{a_1\ast\cdots\ast a_n}, R_{a_1\ast\cdots \ast a_n})\circ
(R_a,  \ldots, R_a)\circ (\phi,  \ldots, \phi),
$$
such that\\
i-\ \ \ \ \ $a_1^n\in H$, \\
ii-\ \ \ \ $\phi:(G, \cdot)\to (H, \ast)$ is an ordinary homomorphism,\\
iii-\ \ \ $a\in C_H(\phi(G))$,\\
iv-\ \ \ \ $a^{n-1}=1$.\\
\end{theorem}

{\bf Proof}. Let $T=(\alpha_1, \alpha_2, \ldots, \alpha_{n+1})$ be a homotopy from $der^n(G, \cdot)$ to $der^n(H, \ast)$. So, we
have
$$
\alpha_{n+1}(x_1x_2\ldots x_n)=\alpha_1(x_1)\ast \cdots \ast \alpha_n(x_n).
$$
Hence, for an arbitrary $x$ and $i$, we have
\begin{eqnarray*}
\alpha_{n+1}(x)&=&\alpha_{n+1}(1^{i-1}\cdot x\cdot 1^{n-i})\\
               &=&\alpha_1(1)\ast \cdots \ast\alpha_{i-1}(1)\ast\alpha_i(x) \ast \alpha_{i+1}(1)\ast \cdots \ast \alpha_n(1).
\end{eqnarray*}
Therefore,
$$
\alpha_i(x)=(\alpha_1(1)\ast \cdots \ast \alpha_{i-1}(1))^{-1}\alpha_{n+1}(x)(\alpha_{i+1}(1)\ast \cdots \ast \alpha_n(1))^{-1}.
$$
Now, suppose $\alpha_i(1)=a_i$ and $a_1\ast a_2\ast \cdots \ast a_n=d$. We have
\begin{eqnarray*}
\alpha_{n+1}(x_1x_2\ldots x_n)&=&\alpha_1(x_1)\ast \alpha_2(x_2)\ast \cdots\ast \alpha_n(x_n)\\
                              &=&\alpha_{n+1}(x_1)\ast d^{-1}\ast \alpha_{n+1}(x_2)\ast d^{-1}\ast\cdots \\
                              & &\ast \alpha_{n+1}(x_{n-1})\ast d^{-1} \ast \alpha_{n+1}(x_n).
\end{eqnarray*}
Hence, we have
\begin{eqnarray*}
\alpha_{n+1}(x_1x_2\ldots x_n)\ast d^{-1}&=&\alpha_{n+1}(x_1)\ast d^{-1}\ast \alpha_{n+1}(x_2)\ast d^{-1}\ast\cdots \\
                                         & &\ast\alpha_{n+1}(x_{n-1})\ast d^{-1}\ast \alpha_{n+1}(x_n)\ast d^{-1}.
\end{eqnarray*}
Therefore, if we let $\psi=R_d^{-1}\alpha_{n+1}=R_a\alpha_{n+1}$, then $\psi:der^n(G, \cdot)\to der^n(H, \ast)$ is a homomorphism.
Now, for any $i$,
\begin{eqnarray*}
\alpha_i(x)&=&(\alpha_1(1)\ast \cdots \ast \alpha_{i-1}(1))^{-1}\alpha_{n+1}(x)(\alpha_{i+1}(1)\ast \cdots \ast \alpha_n(1))^{-1}\\
           &=&(a_1\ast \cdots \ast a_{i-1})^{-1}\ast \psi(x)\ast d\ast (a_{i+1}\ast\cdots\ast a_n)^{-1}\\
           &=&(a_1\ast \cdots \ast a_{i-1})^{-1}\ast \psi(x)\ast (a_1\ast \cdots\ast a_i)\\
           &=&L_{(a_1\ast\cdots \ast a_{i-1})^{-1}}R_{a_1\ast\cdots \ast a_i}\psi(x)\\
           &=&L_{a_i}I_{a_1\ast \cdots \ast a_i}\psi(x).
\end{eqnarray*}
Now, applying lemma 3.1, we obtain the required decomposition for $T$. \\

\begin{corollary}
Suppose $(G, f)=der_{\theta, b}(G, \cdot)$ and $(H, h)=der_{\eta, c}(H, \ast)$ are two $n$-ary groups. Then every homotopy from $(G, f)$
to $(H, h)$ can be decomposed as the  composition of the following homotopies \\

$ (\varepsilon, \theta, \theta^2, \ldots, \theta^{n-1}, R^{-1}_b)$\\

$(\phi, \phi, \ldots, \phi)$\\

$(R_a, R_a, \ldots, R_a)$\\

$(L_{a_1}I_{a_1}, L_{a_2}I_{a_1\ast a_2}, \ldots, L_{a_n}I_{a_1\ast\cdots\ast a_n}, R_{a_1\ast\cdots \ast a_n})$\\

$(\varepsilon, \eta^{-1}, \eta^{-2}, \ldots, \eta^{-{n-1}}, R_c)$\\
such that\\
i-\ \ \ \ \ $a_1^n\in H$, \\
ii-\ \ \ \ $\phi:(G, \cdot)\to (H, \ast)$ is an ordinary homomorphism,\\
iii-\ \ \ $a\in C_H(\phi(G))$,\\
iv-\ \ \ \ $a^{n-1}=1$.\\
\end{corollary}

{\bf Proof}. It is clear that $T=(\varepsilon, \theta, \theta^2, \ldots, \theta^{n-1}, R^{-1}_b)$ is an isotopy from $(G, f)$ to $der^n(G, \cdot)$ and
also $(\varepsilon, \eta^{-1}, \eta^{-2}, \ldots, \eta^{-{n-1}}, R_c)$ is an isotopy from $der^n(H, \ast)$ to $(H, h)$. Now, using 3.3,
the result follows.\\

\begin{corollary}
Let $\psi: der_{\theta, b}(G, \cdot)\to der_{\eta, c}(H, \ast)$ be a homomorphism. Then there exists $a\in H$ and an ordinary homomorphism $\phi:
(G, \cdot)\to (H, \ast)$, such that $\psi=R_a\phi$.
\end{corollary}

{\bf Proof}. The required decomposition for $\psi$ can be obtained by the composition of the first entries of the five
homotopies of the previous corollary.\\

Although there exist some relations between $a$ and $\phi$ in the previous corollary, we are not able to determine these relations completely. However,
in the special case, when $\eta=\varepsilon$ and $c=1$, we have a sufficient and necessary condition for $a$ and $\phi$.

\begin{theorem}
Let $\psi: der_{\theta, b}(G, \cdot)\to der^n(H, \ast)$ be a homomorphism. Then there exists $a\in H$ and an ordinary homomorphism $\phi:
(G, \cdot)\to (H, \ast)$, such that $\psi=R_a\phi$. Further $a$ and $\phi$ satisfy the following conditions;\\
$$
a^{n-1}=\phi(b)\ \ \ and \ \ \ \phi\theta=I_{a^{-1}}\phi.
$$
Conversely, if $a$ and $\phi$ satisfy the above two conditions, then $\psi=R_a\phi$ is a homomorphism $der_{\theta, b}(G, \cdot)\to der^n(H, \ast)$.
\end{theorem}

{\bf Proof}. Let $(G, f)= der_{\theta, b}(G, \cdot)$ and suppose $\psi$ is given. Then by the corollary 3.4, $\psi=R_a\phi$, for some ordinary homomorphism $\phi$ and some $a\in H$.
We have
\begin{eqnarray*}
\psi(b)&=&\psi(1^{n-1}\cdot b)\\
       &=&\psi(1)^n,
\end{eqnarray*}
hence, $a^{n-1}=\phi(b)$. Further,
\begin{eqnarray*}
\phi(\theta(x))\phi(b)&=&\phi(\theta(x)b)\\
                      &=&\psi(\theta(x)b)\ast a^{-1}\\
                      &=&\psi(f(1, x, 1, \dots, 1))\ast a^{-1}\\
                      &=&a\ast \psi(x)\ast a^{n-2}\ast a^{-1}\\
                      &=& a\ast \phi(x)\ast a^{n-2}.
\end{eqnarray*}
So, $\phi\theta=I_{a^{-1}}\phi$. Conversely, suppose $\psi=R_a\phi$, such that $a$ and $\phi$ satisfy the above mentioned
conditions. Now, for all natural number $i$ and all $x$, we have $\phi(\theta^i(x))= a^i\ast \phi(x)\ast a^{-i}$, so
\begin{eqnarray*}
\psi(f(x_1^n))&=&\psi(x_1\theta(x_2)\ldots\theta^{n-1}(x_n)b)\\
              &=&\phi(x_1\theta(x_2)\ldots\theta^{n-1}(x_n)b)\ast a\\
              &=&\phi(x_1)\ast \phi(\theta(x_2))\ast \cdots\ast \phi(\theta^{n-1}(x_n))\ast \phi(b)\ast a\\
              &=&\phi(x_1)\ast a\ast \phi(x_2)\ast a^{-1}\ast a^2\ast \phi(x_3)\ast a^{-2}\ast\cdots \\
              & &\ast a^{n-1}\ast \phi(x_n)\ast a^{-(n-1)}\ast a^{n-1}\ast a\\
              &=&(\phi(x_1)\ast a)\ast(\phi(x_2)\ast a)\ast \cdots\ast (\phi(x_n)\ast a)\\
              &=&\psi(x_1)\ast\cdots\ast \psi(x_n).
\end{eqnarray*}
This completes the proof.\\

There is one more special case which we can determine completely the structure of homomorphisms. In the last theorem of this section, we
investigate this special case.

\begin{theorem}
Let $(G, f)=der_{\theta, b}(G, \cdot)$ and $(H, h)=der_{\eta, c}(H, \ast)$ and suppose further, $(H, \ast)$ is abelian.
Let $\psi: (G, f)\to (H, h)$ be a homomorphism. Then there exists $a\in H$ and an ordinary homomorphism $\phi:
(G, \cdot)\to (H, \ast)$, such that $\psi=R_a\phi$, and
$$
h(\stackrel{(n)}{a})=a\ast \phi(b)\ \ \ and\ \ \ \phi\theta=\eta\phi.
$$
Conversely, if $a\in H$ and $\phi:(G, \cdot)\to (H, \ast)$ is a homomorphism, satisfying the above two conditions, then
$\psi=R_a\phi$ is a homomorphism $(G, f)\to (H, h)$.
\end{theorem}

{\bf Proof}. By the lemma 3.1, and since $(H, \ast)$ is abelian,  every homomorphism $der^n(G, \cdot)\to der^n(H, \ast)$ can be
uniquely decomposed as $R_a\phi$,
such that $a\in H$, $a^{n-1}=1$ and $\phi$ is an ordinary homomorphism. Hence, by 3.3, every homotopy from $(G, f)$ to $(H, h)$
has the form
\newpage
$$
T=(\varepsilon, \eta^{-1}, \ldots, \eta^{-(n-2)}, \varepsilon, R_c)\circ(R_{a_1}, \ldots, R_{a_n}, R_{a_1\ast\cdots\ast a_n})
$$
$$
\circ(\phi, \ldots, \phi)\circ(\varepsilon, \theta, \ldots, \theta^{n-1},R^{-1}_b).
$$
Here, $a_1^n\in H$ and $\phi$ is an ordinary homomorphism. Now, the homomorphisms $\psi: (G, f)\to (H, h)$ are in one-one correspondence
with homotopies which have equal entries. So,  we assume that the entries of $T$ are equal. Let $a=a_1$. It is easy to see that the equality of the
first and the last entries, implies $a_2\ast\cdots\ast a_n\ast c=\phi(b)$. The equality of the first and the second entries, implies
$a_2=\eta(a)$. Similarly, we have $a_{i+1}=\eta^i(a)$. Now, using $a_2\ast\cdots\ast a_n\ast c=\phi(b)$, we conclude
$h(\stackrel{(n)}{a})=a\ast \phi(b)$. Also, the equality of the first and the second entries, implies $\eta^{-1}\phi\theta=\phi$ and more generally,
$\eta^{-i}\phi\theta^i=\phi$. Conversely, if we assume that $a$ and $\phi$ satisfy the required conditions, we can show that all entries of $T$ are equal,
and so $\psi=R_a\phi$ is a homomorphism.\\

\section{ Representations}
Suppose $(G, f)$ is an $n$-ary group. A map $\Lambda:G \rightarrow
GL_m(\mathbb{C})$ with the property
$$
\Lambda(f(x_1^n))=\Lambda(x_1)\Lambda(x_2)\ldots\Lambda(x_n)
$$
is a {\it representation} of $G$. The function
$$
\chi(x)=Tr\ \Lambda(x)
$$
is called the corresponding {\it character} of $\Lambda$. The number
$m$ is the degree of representation. Note that, $\Lambda$ is a
representation of $(G, f)$, iff it is an $n$-ary homomorphism $(G, f)\to
der^n(GL_m(\mathbb{C}))$.

In \cite{Shah},
representation theory of polyadic groups was studied, but representations
dealt with in that paper were considered to have non-empty kernels. The connection of representations of polyadic groups
to representations of their retract investigated in that paper. Further, it is shown that, one can reduce the representation
theory of polyadic groups to the representation theory of ordinary groups (see \cite{Shah}).
In \cite{Shah2}, the second author studied representations of polyadic groups without
any assumptions on kernels, i.e. the representations he  dealt with in \cite{Shah2},
may have empty kernels, as well. It is proved that there is
a one-one correspondence between the set of irreducible
representations of a polyadic group and the  set of irreducible representations of its so called {\it Post covers}. Using
this correspondence, some well-known properties
of irreducible characters of finite groups (such as orthogonality of characters and so on), are generalized to finite polyadic
groups.

In this section, we can use our knowledge about the structure of homomorphisms of polyadic groups, we just obtained in the section 3, to determine
representations of polyadic groups. Applying theorem 3.5, on the case $\Lambda:der_{\theta, b}(G, \cdot)\to
der^n(GL_m(\mathbb{C}))$, we obtain the following theorem.

\begin{theorem}
Let $\Lambda:der_{\theta, b}(G, \cdot)\to der^n(GL_m(\mathbb{C}))$ be a representation. Then we have $\Lambda=R_A\Gamma$, where $A\in GL_m(\mathbb{C})$,
$\Gamma:(G, \cdot)\to GL_m(\mathbb{C})$ is an ordinary representation and we have
$$
A^{n-1}=\Lambda(b)\ \ \  and\ \ \ \Lambda\theta=I^{-1}_{A}\Lambda.
$$
The converse is also true.
\end{theorem}

\end{document}